\documentclass[12pt]{article}
\usepackage[english]{babel}
\usepackage[T1]{fontenc}
\usepackage{graphicx}
\usepackage{amsmath}
\usepackage{amsfonts}
\usepackage{caption}
\usepackage{subcaption}
\usepackage{hyperref}
\usepackage{booktabs}
\usepackage{multirow}
\usepackage{enumitem}
\usepackage{indentfirst}
\usepackage{authblk}
\usepackage{xcolor}
\usepackage{float}

\begin{document}

\title{Learning from user's behaviour of some well-known congested traffic networks}

\date{}

\author[1]{Cardoso, Isolda}
\author[2,3]{Venturato, Lucas}
\author[1]{Walpen, Jorgelina}
\affil[1]{FCEIA, Universidad Nacional de Rosario}
\affil[2]{Departamento de Matem\'atica, Facultad de Ciencias Empresariales, Universidad Austral, Paraguay 1950, S2000FZF, Rosario, Argentina}
\affil[3]{CONICET, Argentina}

%\address{Pellegrini 250, S2000BTP, Rosario, Argentina.}

%\keywords{90B20 TRAFFIC PROBLEMS, 68T20 ARTIFICIAL INTELLIGENCE, 90C33 VARIATIONAL INEQUALITIES}} 
\normalsize

\maketitle

\begin{abstract}
We consider the problem of predicting users' behaviour in a congested traffic network under equilibrium conditions, known as the traffic assignment problem (TAP). We propose and compare 
two machine learning (ML) approaches for this task. The first couples a fully connected neural network with a fixed-point algorithm that enforces Wardrop equilibrium conditions as a post-training correction 
stage. The second explores a path-based neural network architecture inspired by graph neural network (GNN) formulations, which explicitly 
exploits the topological structure of the traffic network. Both schemes are evaluated on several classical benchmark networks.

%\LV{We also explore a path-based neural architecture inspired by message-passing schemes, where layers represent network nodes and connections are constrained by the directed topology, enabling flow propagation across paths while enforcing demand conservation, in contrast to standard GNN formulations.}

\end{abstract}

\section{Introduction}

Intelligent transport systems planning is a research area that offers a variety of problems that are interesting from a mathematical point of view. In the last decades, due to technology evolution, one could argue that such planning might be enhanced and optimized by the use of satellite navigation devices that offer real-time information. However, as expressed in \cite{Morandi} and analysed in \cite{Klein}, the use of such technologies generates the displacement of congestion from one zone to another. It is believed that such effects are a consequence of not knowing the travel choice patterns of users and making no behaviour prediction. Hence, it is still very useful and valuable to study a deterministic and static version of the problem. With this in mind and the fact that information about the transport network and data from user's past experience is available, we propose a different approach. We expect that with the use of data science techniques and machine learning algorithms, travelers' behaviour could be understood and predicted (see \cite{endtoend}).

The problem of computing a traffic assigned flow for a determined network, demand and parameters is known as Traffic Assingment Problem (TAP), and it is modelled according to a behavioural hypothesis over network's users. The difficulties in anticipating effects of road changes, travelers' behaviour, and traffic infrastructure are extremely related to the complexity of the TAP. From a mathematical point of view, the problem is formulated as an equilibrium problem, where equilibrium flows and least costs have to be determined. Such assignment results depend on demands and network's parameters. For each set of demands and parameters, existence and uniqueness results are guaranteed under certain hypothesis (see for example \cite{Patriksson}). So far, effective algorithms have been developed to compute such equilibrium flows and associated least costs. However, such computations are not enough to understand how travelers interact between them and with the traffic network, consequently, predictions are difficult to make. Sensitivity analysis may be useful in this case. Information obtained from evaluating the directions of change that occur in the flows and travel costs when demands and parameters change, can be utilized to optimize changes according to an objective: reducing congestion, minimizing total travel times, etc. Interesting results can be found in \cite{Bell-Iida}, \cite{Choetal}, \cite{Outrata}, \cite{PatRock}, \cite{Pat04}, \cite{Qiuetal}, \cite{Tobin-Friesz} and \cite{Yen}. Most of them provide theoretical results for sensitivity analysis of traffic equilibria: sufficient conditions for the existence of directional derivatives, gradients or subgradients; characterizations for some of them; and computational formulas for their computations, when possible. In particular, in \cite{Pat04} we can find an example that illustrates the fact that the equilibrium solution for the TAP is not differentiable everywhere. The absence of differentiability and convexity (which also occurs) explain the use of tools from variational analysis for the treatment of the problem in most of the above-mentioned papers. This also suggests the use of non conventional machine learning tools. In addition, numerical proposals for the computations of directional derivatives, gradients, or subgradients in general imply the resolution of many network flow problems, which in medium and large scale problems is computationally expensive. 

In recent years, machine learning methods have been proposed as alternatives to classical algorithms for traffic assignment problems. In the static setting, Liu et al. (\cite{endtoend}, \cite{endtoend25}) developed end-to-end frameworks that enforce equilibrium conditions as differentiable layers within the training procedure, while Jungel et al. \cite{WN} propose WardropNet, an equilibrium-augmented architecture that combines classical and learning-based layers. Graph-based approaches have also been explored: Liu and Meidani \cite{LiuGN} and Rahman and Hasan \cite{RH} use graph neural networks to learn the mapping from OD demand to equilibrium flows. ML methods have further been extended to the dynamic setting (Graf et al. \cite{Graf}; Fan et al. \cite{Fan}), though the additional temporal dimension introduces challenges beyond the scope of the present work. 
We propose a decoupled two-stage framework that trades joint optimization for computational tractability, as discussed in Section \ref{introML}; and a message-passing scheme based on GNN formulation.

In order to avoid solving a great number of network assignment problems, we propose an alternative approach for predicting users' behaviour. In contrast to approaches of previous decades, the availability of data may make the difference. Many pairs of data $(x,y)$ where $x$ represents demands and networks parameters and $y$ is the equilibrium solution for the features stored in $x$, can be obtained and stored. If we assume that $y$ is an equilibrium solution that satisfies Wardrop's first principle for a network with characteristics $x$, we are able to design a machine learning tool that benefits from the available $(x,y)$ data to learn users' behaviour and predict future use of the traffic network for new features $x'$. More precisely, we follow the lines of \cite{Wegner}. 

We develop and analyse the effectiveness of a coupled machine learning algorithm that involves a neural network in the first stage and a fixed-point operator in the second. With the neural network we intend, by incorporating the problem restrictions as terms of penalization in the cost function, to approximate the function $ x \overset{\mathcal S}{\rightarrow} y$ which assigns the equilibrium solution for the TAP. The fixed point operator is aimed to adjust the output of the neural network to the equilibrium.

We also explore a neural architecture inspired by message-passing mechanism, where layers represent network nodes and connections are constrained by the directed links, allowing flow propagation along paths while enforcing demand conservation. This design differs from standard graph neural networks, as it does not rely on neighborhood aggregation but instead models directed path-wise flow propagation.

The numerical experiments on some classical networks (such as Small, Steenbrink, Nguyen Dupuis, Sioux Falls, Regular City and Chicago) show that our proposal suitably estimates the solution to the problem of traffic assignment.

In Section \ref{Sec:2}  we establish the mathematical foundations of the TAP. In Section \ref{Sec:3} we give the details of the machine learning approaches. In Section \ref{Sec:4} we give the details for the implementation of the workflow and we establish the metrics for the performance evaluation. In Section \ref{Sec:5} we present the results of the numerical experiments. Finally in Section \ref{Sec:6} we summarize and discuss the work. We also provide additional material in Appendix \ref{App} with supporting figures and tables.

The contribution of this work has the following main aspects: a neural network architecture whose loss function incorporates the feasibility constraints of the TAP; the systematic integration of a post-training fixed-point operator as an equilibrium correction mechanism; the implementation of a message-passing-inspired, path-constrained neural architecture that encodes the directed topology of the network; and finally, numerical tests over several benchmark networks for both proposals.   

%showing that the composition $\mathbb P \circ f$ consistently improves equilibrium conditions with respect to $f$ alone, (where $\mathbb P$ is the fixed point operator and $f$ the neural network predictor) 

\section{The traffic assignment problem}\label{Sec:2}

\subsection{Introduction to the problem}

Let us consider a directed graph $\mathcal G=(\mathcal N, \mathcal A)$, where $\mathcal N$ is the set of nodes and $\mathcal A$ is the set of directed links. Let us also consider the set $\mathcal C\subset \mathcal N\times\mathcal N$ of origin-destination pairs $(p,q)$ and $d_i$ the demand associated to the $ith$ OD pair. The graph represents a transportation network that is assumed to be strongly connected (i.e. there exists at least one route for each OD pair), and each demand $d_i$ has to be assigned to available routes in such a way that every user achieves its objective: arrive from its origin to its destination in minimum time according to Wardrop's first principle (see section \ref{subsec:UE}). In order to obtain a simplified version of the real problem we make the following additional assumptions: travel costs on routes are additive (i.e. the cost of each route is defined as the sum of the costs of the links defining the route); travel costs on links are separable (i.e. the cost of each link is independent of flows on other links); and feasible routes are simple (i.e. no cycles are considered). Finally, we consider link costs functions $t_a$ that are positive, strictly increasing, and parametrized in $m_a, c_a, t^0_a, b_a$. The basic parameters $t^0_a$ and $c_a$ represent the free-flow travel time and the practical capacity of the link $a$, respectively. The positive scalars $m_a$ and $b_a$ are present in the different link performance functions and its interpretation may vary according to the transport network model. The parameters and demands are stored in a feature vector $x$.

We remark that there is an inherent difficulty associated to solving the problem: the equilibrium conditions are formulated in the route flow sense but the feasible set is described in terms of link-flows.

\subsection{  Wardrop's users' equilibrium conditions}\label{subsec:UE}
 
For a deterministic and static version of the problem, the user's equilibrium state is resumed by Wardrop as follows:

\emph{For each OD pair, the journey times on all routes actually used are equal, and less than those which would be experienced by a single vehicle on any unused route}.

The mathematical formulation for the user's equilibrium conditions (DUE, Deterministic User Equilibrium) can be found in \cite{Patriksson}. Briefly, the collection of equalities and inequalities that represent the DUE are listed below. The set of simple feasible routes for each OD pair $(p,q)$ is denoted by $\mathcal R_{pq}$ and $h_r$ is the flow on route $r\in \mathcal R_{pq}$. Then $c_r$ is the travel cost on the route $r$ experienced by each individual user.

\begin{equation}\label{equilconds}
    \begin{array}{rrcll}
     (i) & h_{r}(c_{r}-\pi_{pq})& = & 0, &\forall \;r\in \mathcal R_{pq}, \forall \; (p,q)\in \mathcal C,  \\
     (ii) & c_{r}-\pi_{pq}&\geq &0 , &\forall \;r\in \mathcal R_{pq}, \forall \; (p,q)\in \mathcal C,  \\
     (iii) & \sum \limits_{r\in \mathcal R_{pq}} h_{r} & =& d_{pq}, & \forall \; (p,q)\in \mathcal C,  \\
     (iv) & h_{r}&\geq& 0, & \forall \;r\in \mathcal R_{pq},\forall \; (p,q)\in \mathcal C, \\
     (v) & \pi_{pq}&\geq& 0,& \forall \; (p,q)\in \mathcal C.
    \end{array}
\end{equation}

In this work we consider an equivalent formulation for the DUE as a variational inequality problem. For a set of parameters $m, c, t^0, b \in \mathbb R^{|\mathcal A|}$ and a fixed demand $d \in \mathbb R^{|\mathcal C|}$, an equilibrium flow $h^*$ can be obtained as a solution to the following variational inequality:
\begin{equation}\label{IV}
c(h)^T(h-h^*)\geq 0, \;\; \forall\; h\in H,    
\end{equation}
where $H$ is the product set for $H_{pq}, (p,q)\in \mathcal C$, and $H_{pq}\subset \mathbb R^{|\mathcal R_{pq}|}$ is the set of feasible route flow vectors $h_{r}$, for the OD pair $(p,q)$. The equivalence results are based on the fact that the optimality conditions for the linear program associated with (\ref{IV}) are the conditions that define the DUE, which can also be found in \cite{Patriksson}.

The formulation of the DUE as a variational inequality is the key to design the fixed point operator involved in the second stage of our algorithm.

\subsection{State of the art on proposals for equilibrium flows computation} 

There are several formulations of the DUE as an optimization problem with (\ref{equilconds}) as its optimality conditions. For example a link-route formulation due to Dafermos and Sparrow is given in \cite{DafermosSparrow}, and a link-node formulation due to Beckmann is given in  \cite{Beckmann}. Let us remark that even though both problems are not equivalent since their feasible sets are not the same, their set of equilibria coincide (see \cite{Patriksson}). Moreover, the existing procedures for computing equilibrium solutions are based on these formulations. Frank Wolfe algorithm (FW) and the Dissagregated Simplicial Decomposition (DSD) algorithm are some examples. Both of them are widely used, despite having some drawbacks: the Frank-Wolfe algorithm has as a poor performace for large scale problems, and the DSD needs a great amount of storage capacity for its computations. As we are also concerned with the sensibility analysis for the TAP, we observe that these approaches do not manage efficiently all of their auxiliary computations. Indeed, they fail to take advantage from the vastly computed data.

In this work we consider the Dafermos and Sparrow approach. The resulting optimization problem equivalent to (\ref{equilconds}) is convex as it involves minimizing an integral of an increasing function over a compact domain. 

 Our objective is to approximate the function $\mathcal{S}$ given by $x\in X \to \mathcal{S}(x)=\text{argmin}(\text{TAP})$, in order to be able to make predictions without solving repeatedly new TAPs. Although $\mathcal{S}$ is not differentiable and is not convex (see \cite{Pat04}), it is continuous with respect to the network parameters (see \cite{DafermosNagurney}). Consequently, $\mathcal{S}$ can be approximated by some neural network according to the Universal Approximation Theorem.

\section{Machine learning approach}\label{Sec:3}

\subsection{Introduction}\label{introML}

A transport network $\mathcal G$ with features $\rho$ and equilibrium state $\mathcal{S}(\rho)$ according to the notation in \cite{Pat04}, is a problem that can be revisited from a machine learning point of view, when considerable amount of data is available. According to Wegner \cite{Wegner}, three perspectives can be considered: approximation perspective, probability perspective, and optimization perspective. In this work we propose a combination of the first and the third one in order to accomplish the task of assigning a function $f:X\rightarrow Y$ to a dataset. We expect that $f$ approximates well the function $\mathcal{S}$ in the sense that $f(x) \approx \mathcal{S}(x)$ for every $x\in X$. Nonetheless, we cannot guarantee that the equilibrium conditions are satisfied by the predicted output. The fixed point operator will contribute to this final objective.

Unlike Liu et al. \cite{endtoend}, who enforce equilibrium conditions during training, we decouple the two stages, trading end-to-end optimality for computational tractability.
Equilibrium conditions are very restrictive, making them part of the learning process represented computational effort that could be avoided with the decoupling scheme. Technical details are given in Section \ref{Workflow} and satisfying results are shown in Section \ref{Sec:5}.
Specifically, we propose a ML pipeline to learn feasible link flows from data: existing UE link flows for different demands and network characteristics; and a fixed point operator to fit to equilibrium the learned flows. 

Something similar occurs when we compare with WardropNet approach \cite{WN}. The authors' objective is to embed the equilibrium problem as a final layer into the neural network architecture using an end-to-end approach also. In this case the neural network is used to learn the  parametrization
of latency functions, different from our problem where those parameters are part of the available information. 

Our second proposal, inspired by GNN, is aimed to exploit the traffic network topology avoiding the explicit computation of the paths, hence, this self-supervised scheme is closely related to the nature of the traffic system.

%\LV{We also explore a path-based neural architecture inspired by message-passing schemes, where layers represent network nodes and connections are constrained by the directed topology, enabling flow propagation across paths while enforcing demand conservation, in contrast to standard GNN formulations.}

\subsection{Neural networks design and implementation}

In this section, we present two machine learning models designed to predict unknown flows in a network $\mathcal G$.

\subsubsection{Two-stage Model}

The dataset is constructed as follows: we generate a synthetic dataset by uniformly sampling a set of demands and computing the corresponding arc flow solutions using the Frank-Wolfe algorithm.
This dataset is denoted as $(X,Y)$, where $X$ contains only the random demands $x$. In this study, we consider a simplified version of the traffic assignment problem, fixing specific parameters associated with the traffic network example. We emphasize that the neural network outputs represent route flows, enabling the possibility of testing the equilibrium conditions \eqref{equilconds}.

Considering realistic user behaviour and the exponential increase in computational cost, we assume that excessively long routes are not utilized. Consequently, we select a subset of routes $R_{pq}\subset \mathcal R_{pq}$ in the neural network design. For each network example such subset results from choosing a bound for the maximum length of all the routes.

More precisely, we denote by $K$ such bound. Then $|r|\leq K$, where $|r|$ amounts for the total number of arcs that define the route $r$. 

We build an OD-route matrix $\Delta_{OD}:=[\delta_{pq,r}]$ and an arc route matrix $\Delta_{arc}:=[\delta_{a,r}]$ where $$\delta_{pq,r}=\begin{cases}
1 & \text{ if $r\in R_{pq}$,}\\
0 & \text{ otherwise,}
\end{cases} \; \text{ and } \; \delta_{a,r}=\begin{cases}
1 & \text{ if $a\in r$,}\\
0 & \text{ otherwise.}
\end{cases}$$

We propose a fully connected neural network with an input layer that takes data $x\in X$ and an output layer that predicts route flow vectors $h$. Additionally, we define a custom loss function as follows:
\begin{equation}\label{lossfc}
\mathcal{L}(x)=\|\Delta_{OD}h-x\|^2+\|\Delta_{arc}h-y\|^2,
\end{equation}
where $y\in Y$ represents the arc flow solution for $x$ obtained via the Frank-Wolfe algorithm \cite{Patriksson}. The first term of the loss function measures the feasibility of the predicted route flows, while the second term quantifies the error between the arc flows corresponding to the predicted route flows and the true arc flow solutions.

We remark that, based on the loss function, equilibrium conditions \eqref{equilconds} are not necessarily guaranteed. To address this, we apply a fixed-point operator, proposed in \cite{endtoend}, to the neural network's output, refining the predictions to better satisfy these conditions. The adjusted predictions are then evaluated to verify whether they fulfill the equilibrium conditions. 
For the reader's benefit, we provide a schematic outline of the fixed point operator.

We decompose $H$, the set of all feasible route flow vectors $h$ for all OD pairs $(p,q)$, into two sets defined as follows:
$$H_1:=\{h\in H|h\succeq 0\},\quad \quad H_2:=\{h\in H| \Delta_{OD}h=x\},$$
where in $H_2$, $x$ is the demand vector associated with $h$. Then, projection operators onto $H_1$ and $H_2$ are defined by:
$$\mathbb{P}_{H_1}(h):=h_+, \quad \quad \mathbb{P}_{H_2}(h):=h - \Delta_{OD}^+(\Delta_{OD}h-x),$$
where $\Delta_{OD}^+$ is the Moore–Penrose inverse of $\Delta_{OD}$.
We consider the following fixed point operator in $H$
\begin{equation}\label{fixedpoint}
\mathbb{P}(h):=h-\mathbb{P}_{H_1}(h)+\mathbb{P}_{H_2}(h)\left(2\mathbb{P}_{H_1}(h)-h-\alpha c(\mathbb{P}_{H_1}(h))\right),
\end{equation}
where $c(\cdot)$ is the travel cost function and $\alpha$ is a positive constant.

In particular, if $c(\cdot)$ is maximal monotone, by \cite{Heaton}, the variational inequality \eqref{IV} can be equivalently formulated in terms of the auxiliary fixed point problem $h=\mathbb{P}(h)$. More precisely, $h$ is a solution to \eqref{IV} if and only if $h=\mathbb{P}(h)=\mathbb{P}_{H_1}(h)$.

\subsubsection{Self-supervised Model}\label{MP}

For the message-passing-inspired architecture, we re-arrange the demand information $X$ in the following way: we consider the dataset $(X_{ext}^O,X_{ext}^D)\in \mathbb R^{|\mathcal N|}\times \mathbb R^{|\mathcal N|}$, where $X_{ext}^O$ storages the number of trip intentions from each origin and $X_{ext}^D$ does the same for each destination.

The neural network architecture is defined over successive layers that represent the nodes of the network, with connections constrained by the directed links of the graph. In this setting, flow conservation and propagation are enforced implicitly by the architecture and the loss function, and no explicit OD-route or arc-route incidence matrices are required.
Moreover, the training procedure does not need supervision from equilibrium flows which may not be available. Instead, the loss function penalises deviations from demand satisfaction at destination nodes and minimises TAP's objective function cost.

Finally, the depth of the network is chosen in order to guarantee that the number of layers is equal to or greater than the number of links of the shortest path of each od pair. This ensures that information can propagate along directed paths.

\section{Implementation details}\label{Sec:4}

%\subsection{Neural network architecture}

In this section we provide the notation and describe the architecture of the models presented previously.

We emphasize that our objective is to approximate the true function $\mathcal{S}$, which returns, for every demand vector $x$, the equilibrium arc flow vector.

To this aim, for each example network $\mathcal{G}=(\mathcal{N},\mathcal{A})$, we define:
\begin{itemize}

\item $n_d=|\mathcal{C}|$, the number of OD pairs.

\item $n_a=|\mathcal{A}|$, the total number of arcs.

\item $n_r=\left|\bigcup\limits_{(p,q)\in \mathcal{C}} R_{pq}\right|$, the number of routes considered ($R_{pq}\subseteq \mathcal{R}_{pq}$, for each $(p,q)\in \mathcal{C}$),

\item $K$ the maximum length for all routes in the network.

\end{itemize}

For each OD pair $(p, q)$, we consider a continuous feature $x_{pq}$, which represents the demand associated with that OD pair. Let us recall that we assumed link-additive travel costs on routes, governed by classical link cost functions defined on each arc. In our case:
\begin{equation}\label{lpf}
t_a(y_a) = t_a^0 \left( 1 + 0.15 \left( \frac{y_a}{c_a} \right)^{m_a} \right).
\end{equation}

Related to the two-stage model, we define a deep fully connected neural network with an input layer of size $n_d$, two hidden layers, a ReLu type layer and a normalization layer, both of size $n_r$, and an output layer also of size $n_r$. The output is intended to represent route flow vectors. The model is trained using Adam-type optimizers and the loss function defined in \eqref{lossfc}.

The neural network is followed by a fixed point operator. That is, if $f$ denotes the neural network predictor and $\mathbb{P}$ the fixed point operator described in \eqref{fixedpoint}, we define the composite function $f_0=\mathbb{P}\circ f$ which serves as an approximation of the true function $\mathcal{S}$. Indeed, it is sufficient to choose the exponent $m_a$ in \eqref{lpf} as an odd positive integer number which ensures that the cost function $c(\cdot)$ is maximal monotone, as needed. Finally, we evaluate the accuracy of the approximation $f_0$ provided by the trained model.  

In the case of the self-supervised model, we consider a neural network that mimics the graph structure. In particular, for each OD pair, we define input and output layers of size $|\mathcal{N}|$, that is, node layers, and a number $l$ of hidden node layers according to the expected path lengths for the OD pair, and the connections corresponding to the existing arcs in the graph. Then, for each OD pair we consider this structure as a sub-neural network, sharing the input layer and the loss function. As mentioned in the previous section, demands $x$ are used to generate the input and output data $(X_{ext}^O,X_{ext}^D)$. The output of the sub-neural networks is thus designed to represent the distribution of demand across nodes for its associated OD pair, and the loss function is computed by aggregating their contributions, where each sub-network provides the link flow disaggregated by demand. This architecture bears a conceptual resemblance to column 
generation algorithms, where feasible routes are not enumerated 
explicitly but rather introduced incrementally as needed. 
Analogously, the depth $l$ of the hidden layers implicitly 
controls the set of paths that can be represented by the network: 
each additional layer allows the flow to traverse one additional 
arc, so that paths of length at most $l$ are implicitly captured 
by the architecture without requiring their explicit enumeration. 
This makes the proposed structure particularly well-suited for 
large-scale networks, where the number of feasible routes grows 
exponentially with network size.

\subsection{Workflow}\label{Workflow}

We consider two phases. We start with generating data using SciLab and follow by model training and evaluation using Python, as described below:
\paragraph{Data Generation: SciLab block.}
To train and evaluate the model, we generate data samples $(x, y, c)$, where $x \in \mathbb{R}^{n_d}$ represents a demand vector and $y \in \mathbb{R}^{n_a}$ and $c\in\mathbb{R}$ the corresponding arc flow vector and TAP's objective function value computed via the Frank-Wolf algorithm in Scilab.

\begin{itemize}
    \item For each traffic network, we generate and store a dataset of $N$ demand vectors $\{x^i\}_{i=1}^N$, by uniformly sampling each OD-pair demand $x_{pq}$  in a predefined interval $[x_{\min}, x_{\max}]$. Additionally, we store relevant parameters associated to the traffic network, that are later used by the Python-based algorithms.
    
    We also construct an additional set of demand vectors by sampling each $x_{pq}$ from an expanded interval $[x_{\min}-j_1,, x_{\max}+j_2]$, in order to evaluate the model under out-of-range (extrapolation) conditions.

    \item For each generated demand vector $x^i$, we compute and store the corresponding arc flow vector $y^i$ and the TAP's objective function value by solving the traffic assignment problem using the Frank-Wolfe algorithm implemented in SciLab, with the cost function described in the previous section.
\end{itemize}

\paragraph{Training and evaluation: Python block.}

\begin{itemize}
    \item We import the previously generated data from Scilab, and construct the datasets $(X, Y)$ and $(X_{ext}^O,X_{ext}^D)$ corresponding to each proposed model, as well as $C$, which contains the TAP's objective function value. For the training process, the dataset is divided into training and test subsets. The proportions of each of this subsets are chosen taking into account size and complexity of the traffic net example. %\LV{EN LA SECCIÓN DE EXPERIMENTOS RETOMAR NOMBRANDO PORCENTAJES ESPECIFICOS.}

    \item The output of the neural network consists of route flow vectors for the two-stage model, and of demand vectors for the message-passing architecture. 
    
    \item We apply the fixed-point operator to the two-stage model outputs to obtain updated flows, which are used as refined predictions of the traffic equilibrium.

    \item Finally, we evaluate the performance of both models using the metrics described in the following subsection, focusing on both the accuracy of the predicted arc flows and the DUE condition given by \eqref{IV}. We also compare our results against linear regression and a neural network trained with an MSE loss.
\end{itemize}

\begin{figure}[h]
    \centering
    \includegraphics[scale=0.35]{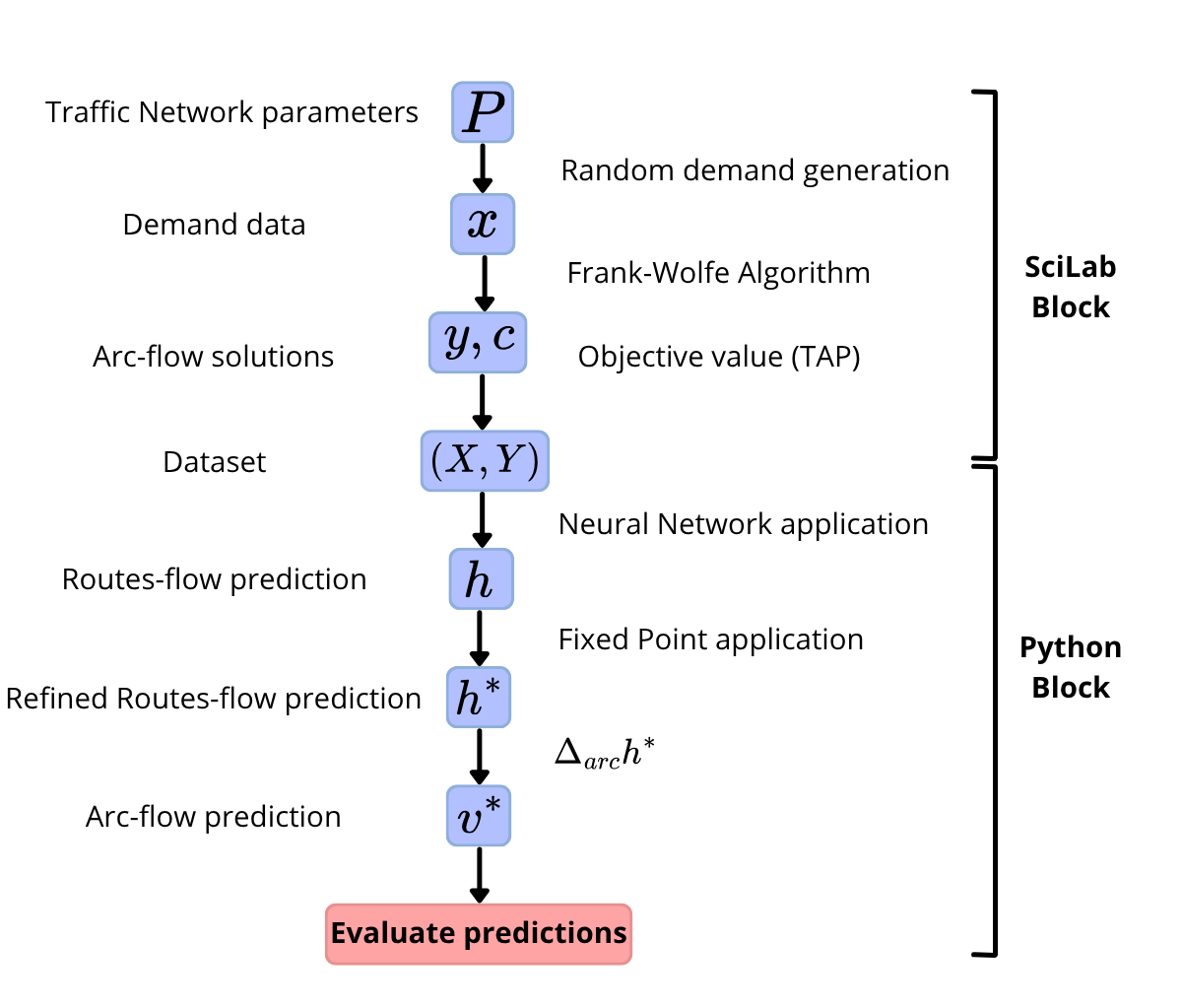}
    \caption{NN+FP Flowchart}
    \label{Flowchart NNFP}
\end{figure}

\begin{figure}[h]
    \centering
    \includegraphics[scale=0.35]{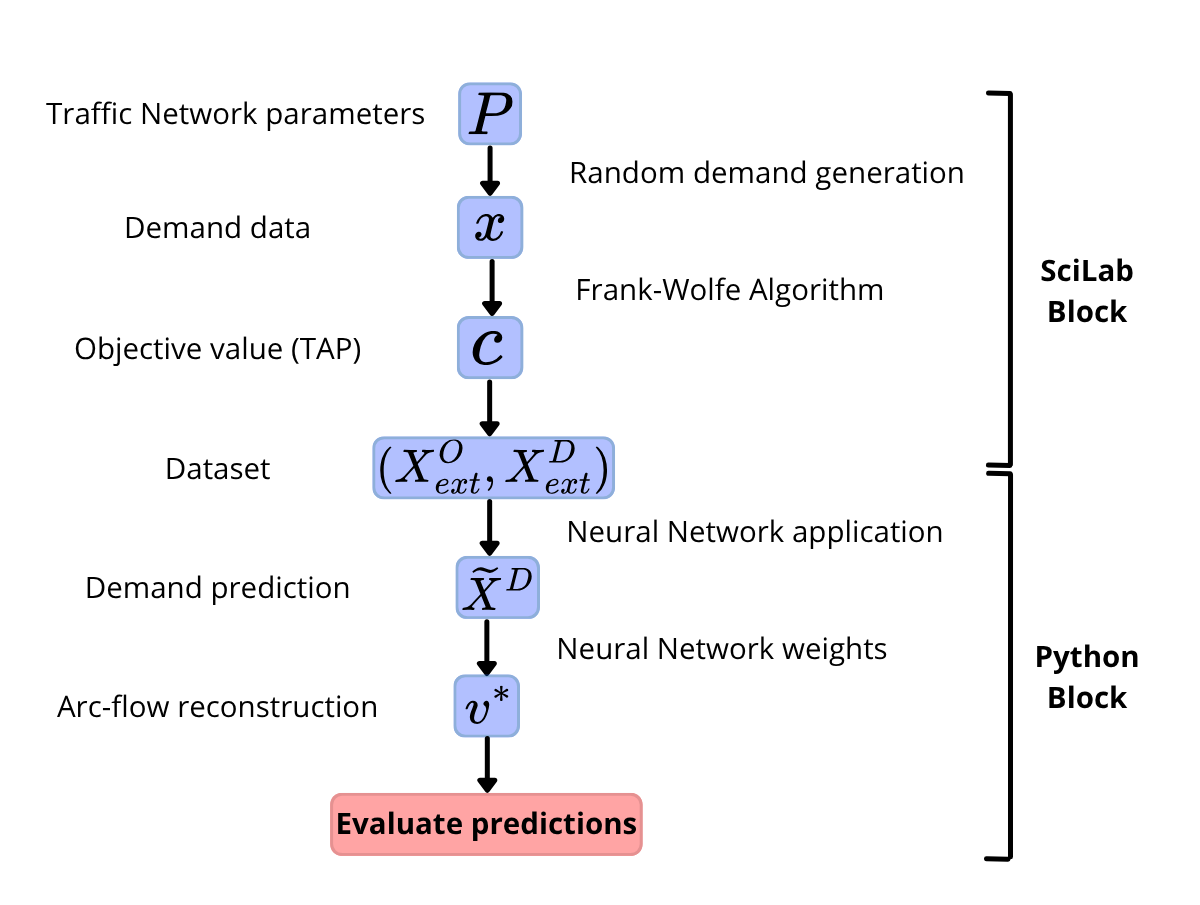}
    \caption{M-P Flowchart}
    \label{Flowchart MP}
\end{figure}

Both complete processes are illustrated in Figures \ref{Flowchart NNFP} and \ref{Flowchart MP}. 

\subsection{Performance}\label{sectionPerformance}

In this section, we present the  error metrics used to evaluate the performance of the both proposed models.
We start by introducing the following notation:

$$pos_0^{ODk} = \textit{ positions where route flow for k-th OD pair equals zero},$$
$$pos_{no}^{ODk} = \textit{ positions where route flow for k-th OD pair is positive},$$
$$c = \textit{ route cost function},$$
and the following conditions
\begin{equation}\label{Mm-cond}
  Mm: \max_i\{c(h)[pos_{no,i}^{ODk}]\}=\min_i\{c(h)[pos_{no,i}^{ODk}]\}, 
\end{equation}
\begin{equation}\label{mM-cond}
  mM: \min_i\{c(h)[pos_{0,i}^{ODk}]\}\geq \max_i \{c(h)[pos_{no,i}^{ODk}]\}.
\end{equation}

We note that, conditions \eqref{Mm-cond} and \eqref{mM-cond} are associated with \eqref{equilconds}-(i) and \eqref{equilconds}-(ii). In fact, if \eqref{Mm-cond} holds for some $(p,q)\in \mathcal{C}$, then $c_{r}-\pi_{pq}=0$. Additionally, condition \eqref{mM-cond} ensures that if $h_{r}=0$, then the associated route cost is higher than the costs of all used routes. That is, condition \eqref{equilconds}-(ii) holds, provided that \eqref{equilconds}-(i) also holds.

For the two-stage model, we rewrite the loss function \eqref{lossfc} as follows
\begin{equation}\label{lossfcRW}
\mathcal{L}(x)=\|\Delta_{OD}h-x\|^2+\|\Delta_{arc}h-y\|^2 = L_1(x)+L_2(x),
\end{equation}
where $h=f(x)$.
Observe that the term $L_1(x)$ is associated with condition~\eqref{equilconds}-(iii), enforcing flow conservation and thus guaranteeing feasibility, which is a hard constraint for any valid solution of the problem. The term $L_2(x)$, on the other hand, measures the discrepancy between the predicted arc flows and the reference solution obtained via the Frank-Wolfe algorithm; this term acts as a soft regularizer and can be weighted or relaxed without compromising feasibility.

For the message-passing approach, the loss function is defined differently, as the model does not directly output route flows.
The loss function consists of two main components. The first term enforces feasibility, as in the previous model, by comparing the prediction with the demand vector, i.e., computing $\|x_{ext}^D-f(x_{ext}^O)\|$, where $f$ denotes the neural network predictor. The second term accounts for the total travel cost induced by the arc flows than can be reconstructed from the learned network weights.
Therefore, the loss function promotes solutions that simultaneously satisfy demand conservation and minimize the associated TAP's objective function value, without relying on explicit supervision from precomputed flow solutions.

In order to analyse both model's performance, we define the following error metrics
$$E_1(x) = \frac{\max\{|\Delta_{arc}h(x)-y(x)|\}}{\max\{y(x)\}},$$
$$E_2(x) = \frac{\max\{|\Delta_{OD}h(x)-x|\}}{\max\{y(x)\}},$$
%$$E_3(x)=\frac{\displaystyle\sum_i|f(x)_i-y(x)_i|}{n_a}$$
$$E_c(x)=\frac{\displaystyle\max_i\{|c(x)_i-cost(x)_i|\}}{\displaystyle\max_i\{c(x)_i,cost(x)_i\}}$$
%$$\bar{y}(x)=\frac{\displaystyle\sum_i|y(x)_i|}{n_a}$$
$$E_{Mm}(x) = \sum_{k}\frac{|\max_i\{c(h(x))[pos_{no,i}^{ODk}]\}-\min_i\{c(h(x))[pos_{no,i}^{ODk}]\}|}{\max_i\{c(h(x))[pos_{no,i}^{ODk}]\}},$$
$$E_{mM}(x) = \sum_k\frac{\left| \max\left\{\max_i\{c(h(x))[pos_{no,i}^{ODk}]\}-\min_i\{c(h(x))[pos_{0,i}^{ODk}]\},0\right\}\right|}{\max\left\{\max_i\{c(h)[pos_{no,i}]\right\},\min_i\{c(h)[pos_{0,i}]\}\}}.$$
$E_1(x)$ measures the error between the arc flows induced by the predicted route flows $h(x)$ and those obtained via the Frank–Wolfe algorithm for the demand vector $x$, in the $\ell_\infty$ sense, i.e., it captures the maximum deviation across arcs. Analogously, $E_2(x)$ and $E_c(x)$ measure the demand feasibility error and the cost error, respectively, both in the $\ell_\infty$ sense. %Finally, $E_3(x)$ corresponds to the average arc flow error, and $\bar{y}(x)$ corresponds to the sample mean of the arc flow vector.

On the other hand, $E_{Mm}(x)$ and $E_{mM}(x)$ are designed to evaluate the DUE conditions.
The average relative errors over all test data are given by
$$E_{1}^X = \frac{1}{|X|}\sum_{x\in X}E_{1}(x),\quad \quad \quad \quad E_{2}^X = \frac{1}{|X|}\sum_{x\in X}E_{2}(x),$$
$$E_{Mm}^X = \frac{1}{n_d|X|}\sum_{x\in X}E_{Mm}(x),\quad \quad \quad \quad E_{mM}^X = \frac{1}{n_d|X|}\sum_{x\in X}E_{mM}(x).$$
%\LV{$$E_3^X = \frac{1}{|X|}\sum_{x\in X}E_3(x),$$
$$E_c^X = \frac{1}{|X|}\sum_x E_c(x),$$
%$$\bar{Y} = \frac{1}{|X|}\sum_{x\in X}\bar{y}(x),$$}

In the case of the two-stage model, we train our neural network using a dynamic loss function: during the first half of the epochs, we apply only $L_{1}$, and during the second half, we apply full $\mathcal{L}$. This improves computational times in the training phase.
Reversing this schedule does not improve neither the predictions nor the times. 

Additionally, we allow for a margin of error, denoted by $\epsilon$, in the conditions $Mm$ and $mM$, in the testing, which correspond to consider the following conditions $Mm'$ and $mM'$
$$mM': \min_i\{c(h)[pos_{0,i}^{ODk}]\}\geq c(h)[pos_{no,j}^{ODk}])\cdot (1-\epsilon), \quad \forall j,$$
$$Mm': \left|\max_i\{c(h)[pos_{no,i}^{ODk}]\}-\min_i\{c(h)[pos_{no,i}^{ODk}]\}\right|\leq \epsilon\max_i\{c(h)[pos_{no,i}^{ODk}]\},$$
for $0<\epsilon<<1$.

The self-supervised model employs its corresponding full loss function throughout the entire training process.

\section{Numerical Results}\label{Sec:5}

In this section, we present the results obtained by applying our proposed models and linear regression to the benchmark networks. Performance is measured using the error metrics described in the previous section.
Each of the numerical experiments was run with the PyTorch library on a AMD Ryzen 5 CPU 5600h with 8 GB of RAM.

It is worth noting that the linear regression models are trained using an MSE loss on arc flows and, as a consequence, their predictions do not guarantee feasibility. Indeed, satisfactory results are obtained only for small traffic networks, which motivates the use of architectures that explicitly respect the underlying problem structure.

Figures \ref{Small} to \ref{Chicago} show each traffic network with their corresponding demands, and Table \ref{Tabla0} presents the associated network parameters, as mentioned before.

\subsection{Development}\label{Develop}

%We performed an initial approach which lead us to the formulation above. 
Our first attempt considered the use of a classical linear regression method in order to build a regressor that predicts the output for new data. We trained the regressor with the dataset $(X,Y)$. This proposal allowed us to obtain good performance in terms of link flow for the predictions. However, since the method was designed for the arc flow variable, we were not able to check performance in terms of feasibility nor equilibrium conditions. Moreover, predictions failed to satisfy feasibility constraints, except for the case of the Small example. To address this issue, we designed a neural network with output size equal to the total number of routes $n_r$, so that the outputs could be directly interpreted as route flow vectors. We considered the loss function \eqref{lossfc}, but the results were still unsatisfactory in terms of equilibrium conditions. Consequently, we considered  adding equilibrium conditions to the loss function, but the computational time was too high. With all this experience in mind, since the predictor already built by the neural network worked well both in terms of feasibility and approximation to the actual arc flow solution, we thought of imposing the equilibrium conditions as a final stage. We applied the fixed-point procedure only after training the neural network using the loss function \eqref{lossfc} with better results. However, a new issue aroused: training required a large number of epochs. We mitigated it with a dynamic loss function, which led to improved results. 

As an alternative, we explored the message-passing approach, which avoids the need for precomputed route flow data $Y$. In this case, arc flows are implicitly determined by the trainable weights of the network, and the loss function is designed to simultaneously enforce demand feasibility and minimise the total travel cost, without explicit supervision.

\vspace{1cm}

\subsection{Computational performance}

In this section we present the test results for each model, evaluating both predictive performance and computational efficiency.
The performance is measured in terms of the error metrics presented in Section \ref{sectionPerformance}.
The results are summarized in Tables \ref{tab:E1}-\ref{tab:time}.

In Tables \ref{tab:E1} and \ref{tab:E2} we can observe the performance of both proposed methods and linear regression methods that were implemented. Table \ref{tab:E1} shows performance of predicted flows as regards closeness to link flows obtained by Frank-Wolfe algorithm. Table \ref{tab:E2} shows performance of predicted flows in terms of feasibility.

According to the results in Table \ref{tab:E1}, both methods perform acceptably in general terms when predicting flows close to the Frank-Wolfe solution. However, as shown in Table \ref{tab:E2}, satisfactory feasibility of the predictions is achieved only by the NN+FP and M-P methods. This highlights the necessity of explicitly incorporating demand satisfaction constraints into the network design itself. Concretely, being close to an equilibrium flow does not necessarily guarantee the feasibility of the prediction. This observation was precisely  the one that motivated the self-supervised formulation proposed in Section \ref{MP}, where feasibility is enforced directly through the loss function rather than as a post-hoc correction.

Tables \ref{tab:Ec}, \ref{tab:EmM} and \ref{tab:EMm} assess the performance of the methods in terms of equilibrium flow prediction. Specifically, Table \ref{tab:Ec} measures the distance to the optimal value of the TAP objective function obtained via the Frank-Wolfe algorithm. In these terms, all methods produce acceptable predictions.
Tables \ref{tab:EmM} and \ref{tab:EMm} analyse only the NN and NN+FP methods, since the performance metrics considered in this case require route flow information, which is not available for the remaining methods. In particular, Table \ref{tab:EmM}, which examines unused routes with lower cost than those actually used (condition \eqref{equilconds}-(ii)), reflects an improvement in the predictions of NN+FP over NN. Regarding Table \ref{tab:EMm}, which compares the travel costs of routes used in the predictions, greater difficulty is observed in guaranteeing cost equality among them (condition\eqref{equilconds}-(i)). Among the possible causes identified, this may be attributed to the limited number of routes considered in the numerical experiments.

Finally, Table \ref{tab:time} compares the running times of the different models used in the numerical experiments against those of the Frank-Wolfe algorithm. In all cases, the machine learning methods achieve running times several orders of magnitude shorter than those of the Frank-Wolfe algorithm, supporting the computational tractability argument motivating the proposed approach.

\section{Conclusions}\label{Sec:6}

Our original objective was to try to solve the classical traffic assignment problem with machine learning methods, given the large amount of available data and the need to reduce computational cost of the existing algorithms.

The Flowcharts given in \ref{Flowchart NNFP} and \ref{Flowchart MP} summarize the main details of our work. We built two predictors. The first one consists of an ad-hoc neural network combined with a fixed point operator and the second one, a message passing type neural network. We tested our predictors on several traffic networks of the literature and obtained good results in all of the terms we considered, namely feasibility, solution approximation, satisfaction of the equilibrium conditions and execution times.

Some special comments can be made towards the NN+FP method. The proposed method not only performs acceptably but also generates useful information as regards route flows, route costs and route composition without extra computational cost, information that is not available when FW algorithm is used, for example.

As regards the M-P method, it also performs acceptably even though the arc flow data from the Frank-Wolfe algorithm is not needed, and moreover, there is no need of explicit computations of the paths.

Finally, we remark that the NN+FP method is more costly in terms of time than the M-P method.

% \begin{figure}[h!]
%     \centering
%     \includegraphics[scale=0.25]{Flowchart2.png}
%     \caption{Flowchart}
%     \label{Flowchart2}
% \end{figure}

\section*{Acknowledgement}
This work was partially supported by Universidad Nacional de Rosario under Grant No. 80020230300094UR, ``Matemática, aprendizaje automático y aplicaciones al transporte''. Additionally, the second author was partly supported by Universidad Austral under Grant No. 006–25CI2001; and CONICET under Grant No. 11220220100532.

\section*{Disclosure Statement}

The authors report there are no competing interests to declare.

\newpage
\appendix

\section{Appendix: Figures and Tables}\label{App}

\begin{figure}[H]
    \centering
    \includegraphics[scale=0.15]{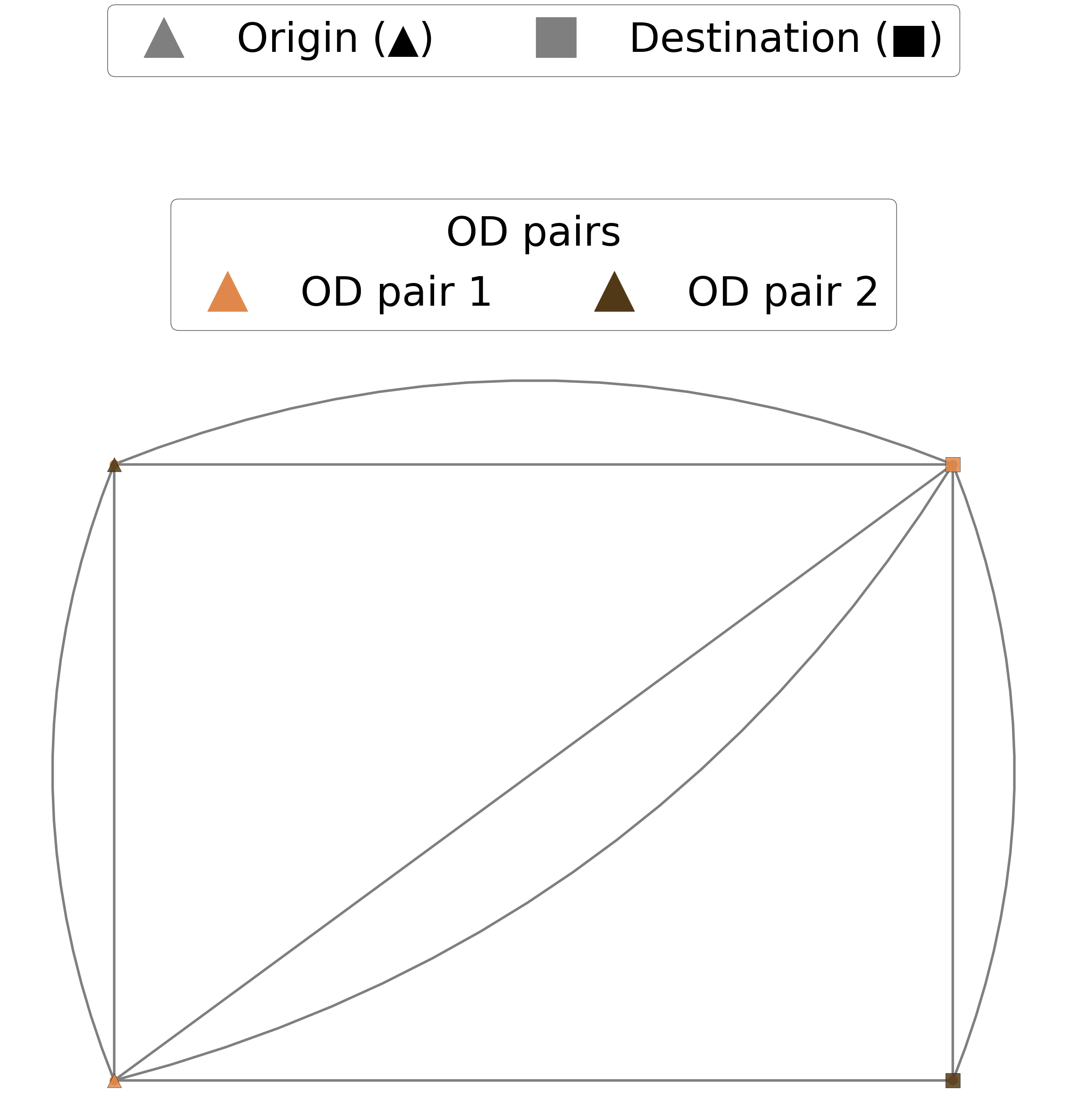}
    \caption{Small}
    \label{Small}
\end{figure}

\begin{figure}[H]
    \centering
    \includegraphics[scale=0.15]{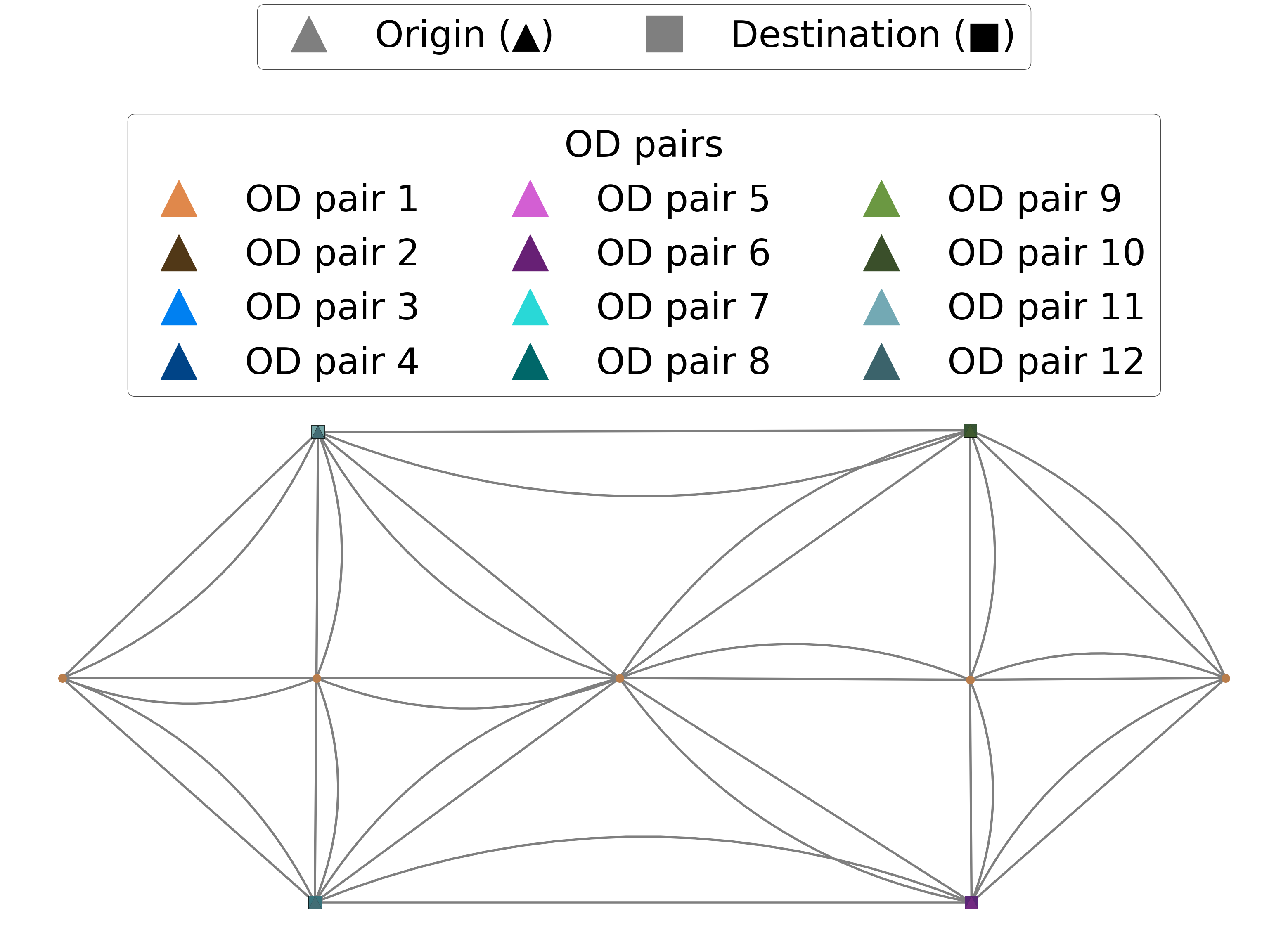}
    \caption{Steenbrink}
    \label{Steenbrink}
\end{figure}

\begin{figure}[H]
    \centering
    \includegraphics[scale=0.15]{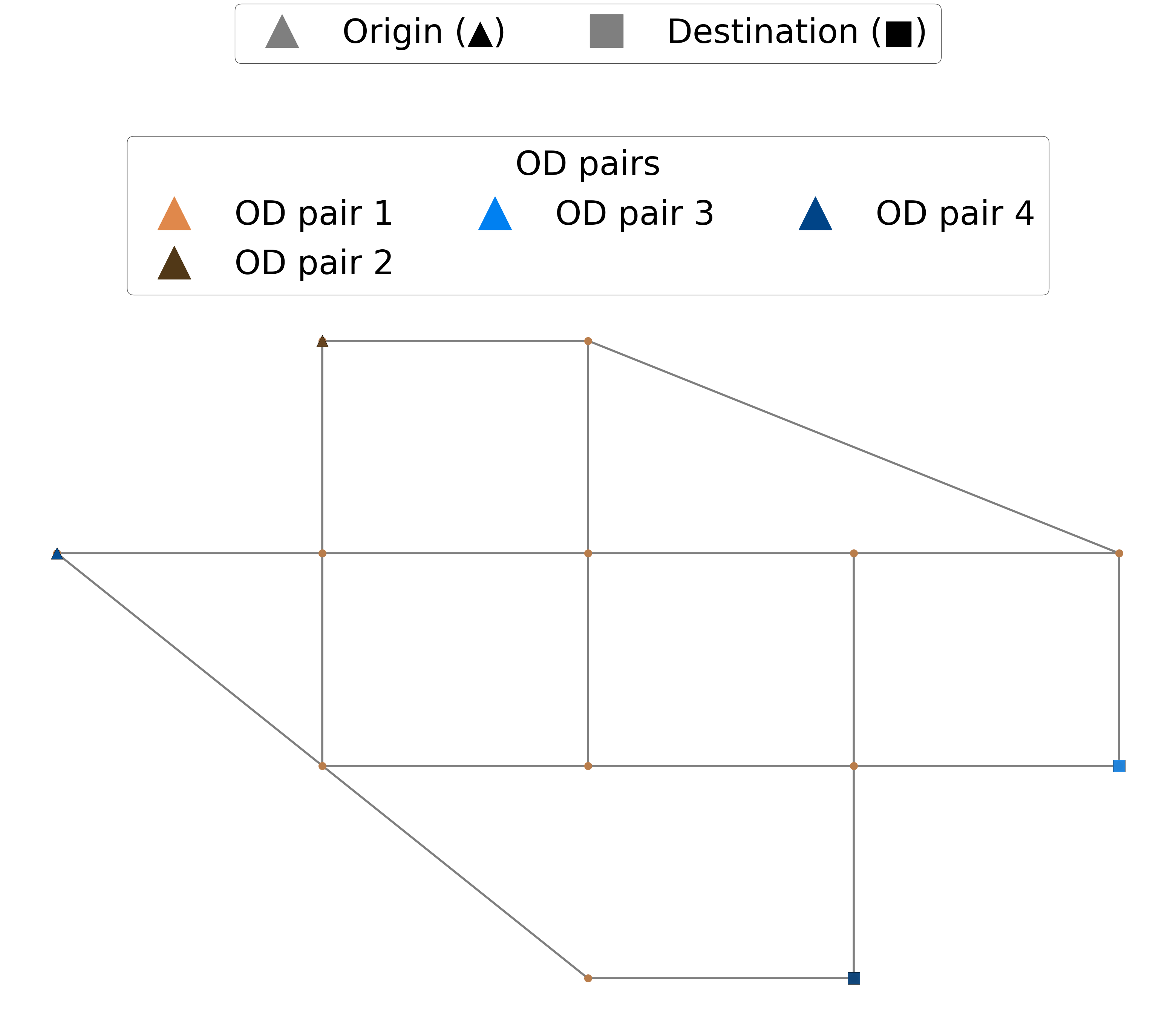}
   \caption{Nguyen Dupuis}
    \label{Nguyen Dupuis}
\end{figure}

\begin{figure}[H]
    \centering
    \includegraphics[scale=0.15]{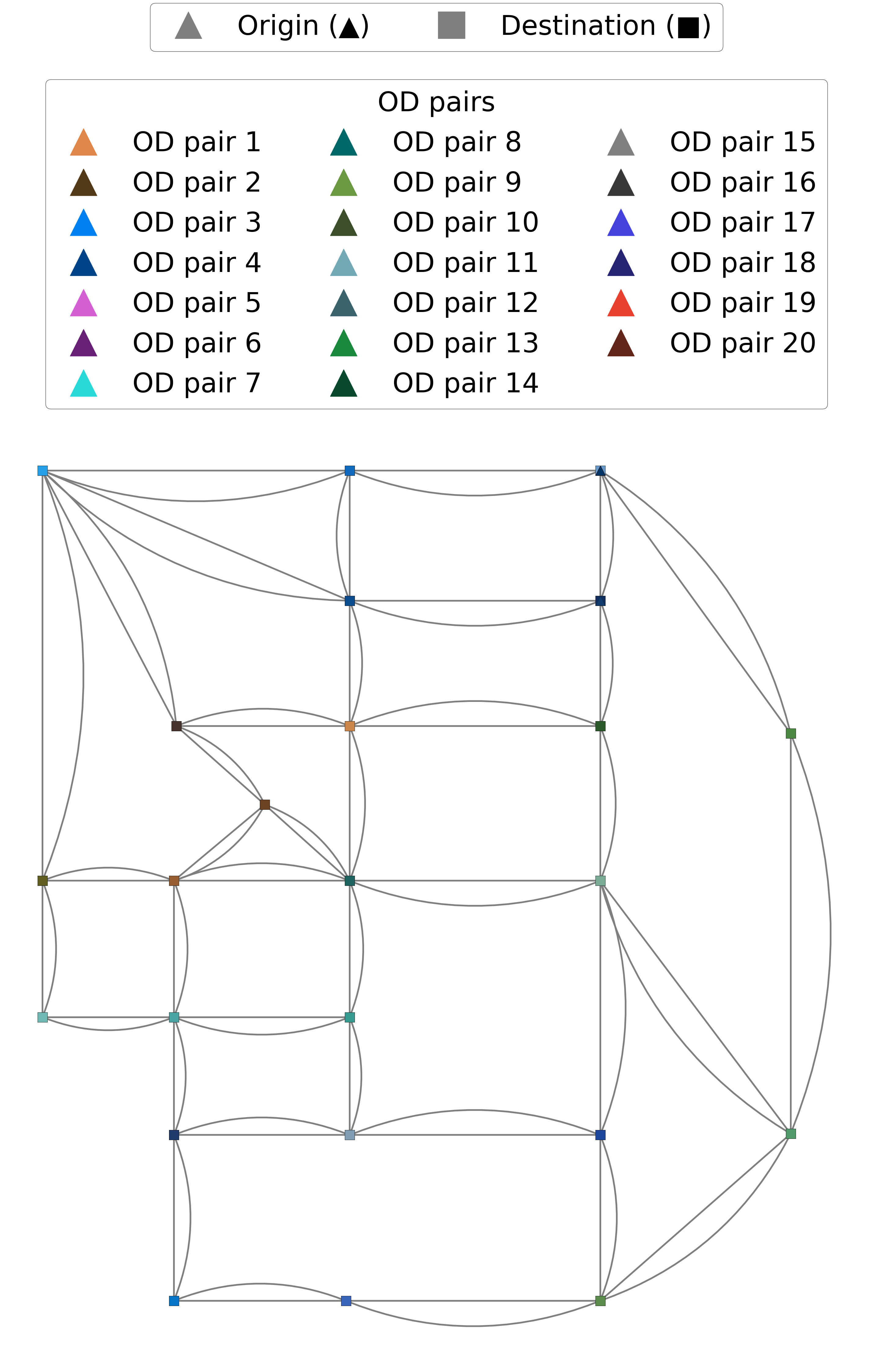}
    \caption{Sioux Falls Example}
    \label{Sioux Falls}
\end{figure}

\begin{figure}[H]
    \centering
    \includegraphics[scale=0.15]{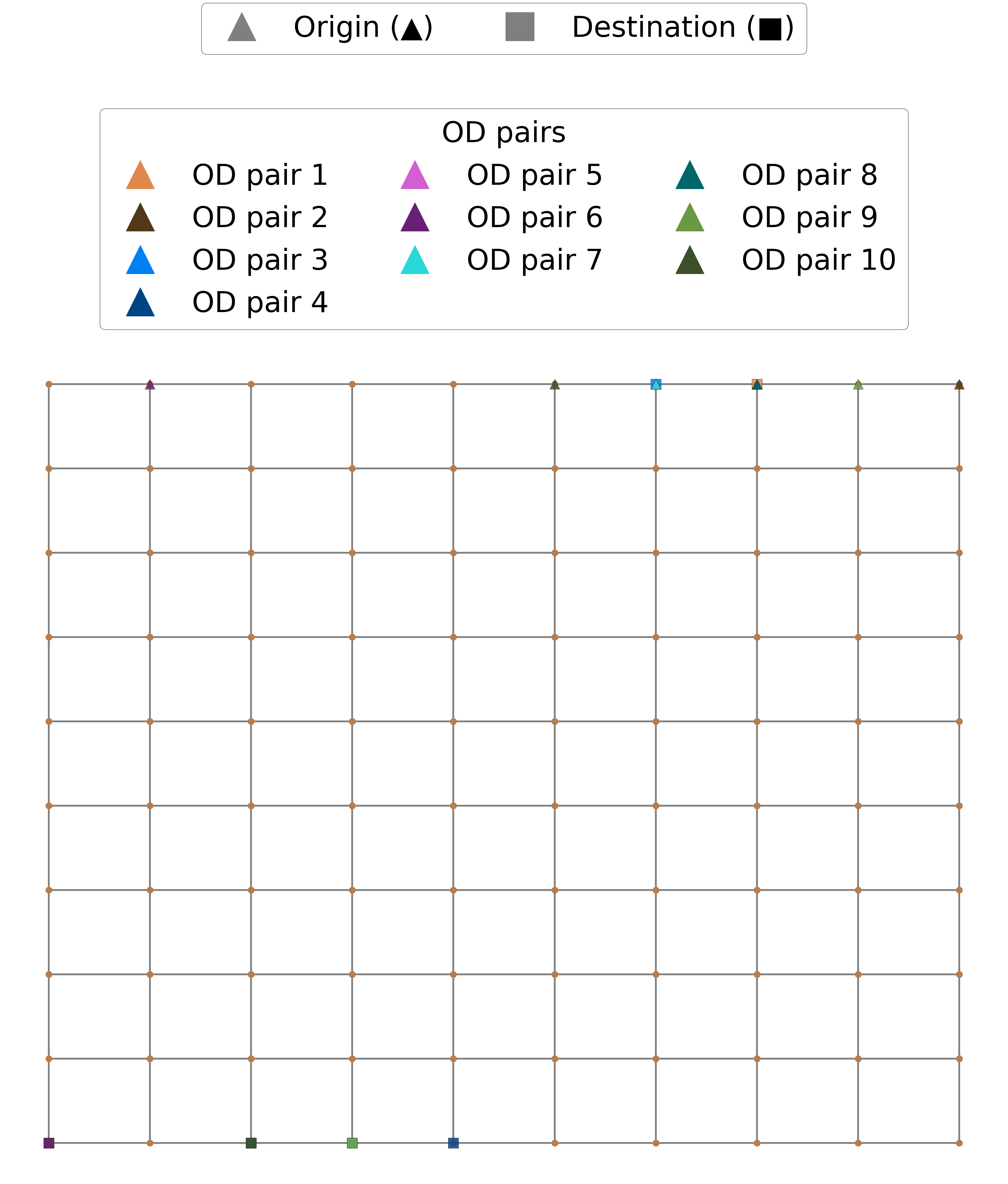}
    \caption{Regular Example}
    \label{Regular}
\end{figure}

\begin{figure}[H]
    \centering
    \includegraphics[scale=0.15]{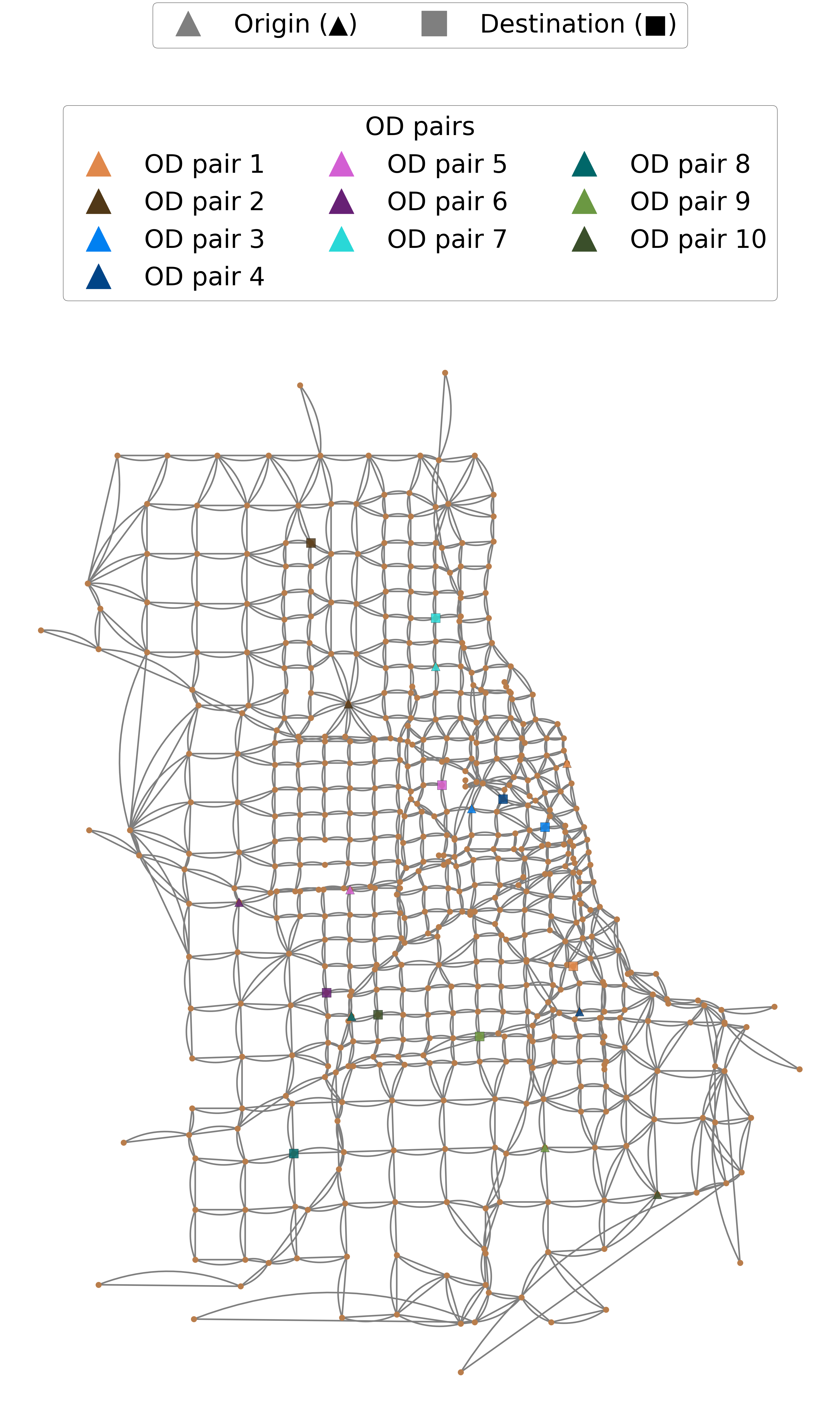}
    \caption{Chicago Example}
    \label{Chicago}
\end{figure}

\newpage

\begin{table}[h!]
\centering
\caption{Traffic networks parameters}
\begin{tabular}{|c|c|c|c|c|c|c|c|}
\hline
& $|\mathcal{N}|$ & $n_d$ & $n_a$ & $n_r$ & $m_{pq}$   \\ 
\hline
Small & 4 & 12 &  9 & 30 & 4 \\ 
\hline
Steen & 9 & 12 & 36 & 1112 & 7  \\ 
\hline
Nguyen & 13 & 4 & 19 & 25 & 5   \\ 
\hline
Sioux & 24 & 100  & 76 & 3541 & 4   \\ 
\hline
Sioux & 24 & 528  & 76 & 21006 & 4   \\ 
\hline
Regular & 100 & 30  & 180 & 3512 & 4  \\ 
\hline
Chicago & 546 & 10 & 2176 & 73520 & 4  \\ 
\hline
\end{tabular}
\label{Tabla0}
\end{table}

\vspace{1cm}

% -------------------------------------------------------
% E_1^X
% -------------------------------------------------------
\begin{table}[h!]
\small
\centering
\caption{Algorithm Performance -- $E_1^X$}
\begin{tabular}{| l | c | c | c | c | c | c | c |}
\hline
\multirow{2}{*}{\textbf{Network}} 
  & \multicolumn{3}{c|}{\textbf{Linear Regression}} 
  & \multirow{2}{*}{\textbf{NN (MSE)}} 
  & \multirow{2}{*}{\textbf{NN}} 
  & \multirow{2}{*}{\textbf{NN+FP}} 
  & \multirow{2}{*}{\textbf{M-P}} \\
\cline{2-4}
& \textbf{Classic} & \textbf{Ridge} & \textbf{Lasso} & & & & \\
\hline
Small                & 0.032 & 0.032 & 0.032 & 0.034 & 0.02  & 0.06  & 0.144 \\
\hline
Small-out            & 0.068 & 0.068 & 0.068 & 0.281 & 0.11  & 0.08  & 0.229 \\
\hline
Steenbrink           & 0.063 & 0.063 & 0.063 & 0.096 & 0.07  & 0.12  & 0.256 \\
\hline
Steenbrink-out       & 0.148 & 0.148 & 0.148 & 0.279 & 0.20  & 0.17  & 0.423 \\
\hline
Nguyen Dupuis        & 0.011 & 0.011 & 0.011 & 0.046 & 0.15  & 0.14  & 0.041 \\
\hline
Nguyen Dupuis-out    & 0.088 & 0.088 & 0.088 & 0.274 & 0.21  & 0.15  & 0.078 \\
\hline
Sioux Falls          & 0.234 & 0.043 & 0.043 & 0.057 & 0.05  & 0.11  & 0.242 \\
\hline
Sioux Falls-out      & 0.713 & 0.28  & 0.095 & 0.112 & 0.15  & 0.27  & 0.312 \\
\hline
Regular 10-10-10     & 0.032 & 0.032 & 0.032 & 0.081 & 0.04  & 0.04  & 0.098 \\
\hline
Regular 10-10-10-out & 0.062 & 0.062 & 0.062 & 0.093 & 0.15  & 0.06  & 0.205 \\
\hline
Regular 10-10-30     & 0.030 & 0.030 & 0.030 & 0.093 & 0.04  & 0.06  & 0.125 \\
\hline
Regular 10-10-30-out & 0.055 & 0.055 & 0.055 & 0.250 & 0.12  & 0.12  & 0.235 \\
\hline
Chicago 10           & 0.0   & 0.0   & 0.0   & 0.069 & 0.11  & 0.10  & 0.576 \\
\hline
Chicago 10-out       & 0.0   & 0.0   & 0.0   & 0.429 & 0.57  & 0.24  & 0.699 \\
\hline
% Chicago 1000         & 0.434 & 0.024 & 0.019 & 0.024 & -- & -- & -- \\
% \hline
% Chicago 1000-out     & 1.121 & 0.057 & 0.046 & 0.043 & -- & -- & -- \\
% \hline
\end{tabular}
\label{tab:E1}
\end{table}

% -------------------------------------------------------
% E_2^X
% -------------------------------------------------------
\begin{table}[h!]
\small
\centering
\caption{Algorithm Performance -- $E_2^X$}
\begin{tabular}{| l | c | c | c | c | c | c | c |}
\hline
\multirow{2}{*}{\textbf{Network}} 
  & \multicolumn{3}{c|}{\textbf{Linear Regression}} 
  & \multirow{2}{*}{\textbf{NN (MSE)}} 
  & \multirow{2}{*}{\textbf{NN}} 
  & \multirow{2}{*}{\textbf{NN+FP}} 
  & \multirow{2}{*}{\textbf{M-P}} \\
\cline{2-4}
& \textbf{Classic} & \textbf{Ridge} & \textbf{Lasso} & & & & \\
\hline
Small                & 0.0 & 0.0 & 0.0 & 0.118 & 0.0  & 0.0   & 0.005 \\
\hline
Small-out            & 0.0 & 0.0 & 0.0 & 0.185 & 0.05  & 0.0   & 0.005 \\
\hline
Steenbrink           & 0.422 & 0.422 & 0.421 & 0.333 & 0.02  & 0.0   & 0.0   \\
\hline
Steenbrink-out       & 0.846 & 0.846 & 0.844 & 0.522 & 0.17  & 0.0   & 0.0   \\
\hline
Nguyen Dupuis        & 0.976 & 0.976 & 0.976 & 0.917 & 0.33  & 0.0   & 0.0   \\
\hline
Nguyen Dupuis-out    & 0.922 & 0.922 & 0.922 & 1.06 & 0.34  & 0.0   & 0.0   \\
\hline
Sioux Falls          & 12.714 & 1.917 & 1.184 & 0.434 & 0.28  & 0.0   & 0.0   \\
\hline
Sioux Falls-out      & 29.021 & 3.842 & 2.292 & 0.804 & 0.44  & 0.12  & 0.0   \\
\hline
Regular 10-10-10     & 1.581 & 1.581 & 1.581 & 1.535 & 0.03  & 0.0   & 0.001 \\
\hline
Regular 10-10-10-out & 1.224 & 1.223 & 1.224 & 1.354 & 0.21  & 0.0   & 0.001 \\
\hline
Regular 10-10-30     & 3.22 & 3.218 & 3.218 & 3.108 & 0.10  & 0.0   & 0.006 \\
\hline
Regular 10-10-30-out & 2.463 & 2.459 & 2.458 & 2.317 & 0.21  & 0.02  & 0.006 \\
\hline
Chicago 10           & 0.0 & 0.0 & 0.0 & 0.127 & 0.16  & 0.06  & 0.003 \\
\hline
Chicago 10-out       & 0.0 & 0.0 & 0.0 & 0.437 & 1.01  & 0.18  & 0.003 \\
\hline
% Chicago 1000         & 237.069 & 184.714 & 184.74 & 181.388 & -- & -- & -- \\
% \hline
% Chicago 1000-out     & 276.018 & 125.397 & 123.379 & 122.893 & -- & -- & -- \\
% \hline
\end{tabular}
\label{tab:E2}
\end{table}

\begin{table}[h!]
\small
\centering
\caption{Algorithm Performance -- $E_c$}
\begin{tabular}{| l | c | c | c | c | c | c | c |}
\hline
\multirow{2}{*}{\textbf{Network}} 
  & \multicolumn{3}{c|}{\textbf{Linear Regression}} 
  & \multirow{2}{*}{\textbf{NN (MSE)}} 
  & \multirow{2}{*}{\textbf{NN}} 
  & \multirow{2}{*}{\textbf{NN+FP}} 
  & \multirow{2}{*}{\textbf{M-P}} \\
\cline{2-4}
& \textbf{Classic} & \textbf{Ridge} & \textbf{Lasso} & & & & \\
\hline
Small                & 0.011 & 0.011 & 0.011 & 0.041 & 0.014 & 0.014 & 0.025  \\
\hline
Small-out            & 0.048 & 0.048 & 0.048 & 0.300 & 0.062 & 0.023 & 0.012  \\
\hline
Steenbrink           & 0.011 & 0.011 & 0.011 & 0.124 & 0.050 & 0.115 & 0.070  \\
\hline
Steenbrink-out       & 0.112 & 0.112 & 0.112 & 0.222 & 0.218 & 0.226 & 0.292  \\
\hline
Nguyen Dupuis        & 0.0   & 0.0   & 0.0   & 0.016 & 0.008 & 0.004 & 0.0002 \\
\hline
Nguyen Dupuis-out    & 0.001 & 0.001 & 0.001 & 0.151 & 0.061 & 0.019 & 0.0011 \\
\hline
Sioux Falls          & 0.076 & 0.019 & 0.082 & 0.080 & 0.081 & 0.366 & 0.018  \\
\hline
Sioux Falls-out      & 0.423 & 0.219 & 0.131 & 0.970 & 0.385 & 0.633 & 0.085  \\
\hline
Regular 10-10-10     & 0.001 & 0.001 & 0.001 & 0.083 & 0.010 & 0.013 & 0.013  \\
\hline
Regular 10-10-10-out & 0.005 & 0.006 & 0.005 & 0.157 & 0.133 & 0.076 & 0.017  \\
\hline
Regular 10-10-30     & 0.001 & 0.001 & 0.001 & 0.030 & 0.052 & 0.042 & 0.055  \\
\hline
Regular 10-10-30-out & 0.003 & 0.003 & 0.004 & 0.034 & 0.156 & 0.135 & 0.074  \\
\hline
Chicago 10           & 0.0   & 0.0   & 0.0   & 0.046 & 0.010 & 0.007 & 0.055  \\
\hline
Chicago 10-out       & 0.0   & 0.0   & 0.0   & 0.152 & 0.508 & 0.204 & 0.059  \\
\hline
% Chicago 1000         & 0.012 & 0.007 & 0.007 & 0.039 & --  & --  & -- \\
% \hline
% Chicago 1000-out     & 0.031 & 0.017 & 0.023 & 0.040 & --  & --  & -- \\
% \hline
\end{tabular}
\label{tab:Ec}
\end{table}

% -------------------------------------------------------
% E_mM^X
% -------------------------------------------------------

\begin{table}[h!]
\centering
\caption{Algorithm Performance -- $E_{mM}^X$ ($\epsilon = 0.1/0.05/0.01$)}
\begin{tabular}{| l | c | c |}
\hline
\textbf{Network} & \textbf{NN} & \textbf{NN+FP} \\
\hline
Small                & 0.004/0.005/0.005 & 0.0/0.0/0.0       \\
\hline
Small-out            & 0.002/0.002/0.002 & 0.0/0.0/0.0       \\
\hline
Steenbrink           & 0.336/0.362/0.383 & 0.0/0.001/0.003   \\
\hline
Steenbrink-out       & 0.537/0.557/0.572 & 0.018/0.031/0.045 \\
\hline
Nguyen Dupuis        & 0.005/0.019/0.031 & 0.0/0.0/0.0       \\
\hline
Nguyen Dupuis-out    & 0.006/0.022/0.035 & 0.0/0.0/0.0       \\
\hline
Sioux Falls          & 0.253/0.28/0.301  & 0.0/0.003/0.01    \\
\hline
Sioux Falls-out      & 0.254/0.28/0.3    & 0.001/0.003/0.008 \\
\hline
Regular 10-10-10     & 0.093/0.129/0.159 & 0.0/0.0/0.005     \\
\hline
Regular 10-10-10-out & 0.146/0.183/0.209 & 0.006/0.012/0.025 \\
\hline
Regular 10-10-30     & 0.092/0.126/0.152 & 0.002/0.004/0.01  \\
\hline
Regular 10-10-30-out & 0.091/0.124/0.151 & 0.016/0.025/0.034 \\
\hline
Chicago 10           & 0.29/0.313/0.330  & 0.140/0.154/0.166 \\
\hline
Chicago 10-out       & 0.306/0.328/0.344 & 0.164/0.179/0.192 \\
\hline
\end{tabular}
\label{tab:EmM}
\end{table}

% -------------------------------------------------------
% E_Mm^X
% -------------------------------------------------------

\begin{table}[h!]
\centering
\caption{Algorithm Performance -- $E_{Mm}^X$ ($\epsilon = 0.1/0.05/0.01$)}
\begin{tabular}{| l | c | c | c |}
\hline
\textbf{Network} & \textbf{NN} & \textbf{NN+FP}\\
\hline
Small                & 0.121/0.128/0.131 & 0.116/0.128/0.130 \\
\hline
Small-out            & 0.152/0.158/0.158 & 0.159/0.164/0.164 \\
\hline
Steenbrink           & 0.430/0.431/0.432 & 0.368/0.368/0.368 \\
\hline
Steenbrink-out       & 0.641/0.641/0.641 & 0.484/0.484/0.484 \\
\hline
Nguyen Dupuis        & 0.056/0.065/0.073 & 0.03/0.036/0.036  \\
\hline
Nguyen Dupuis-out    & 0.047/0.062/0.068 & 0.017/0.032/0.034 \\
\hline
Sioux Falls          & 0.594/0.595/0.596 & 0.569/0.057/0.57  \\
\hline
Sioux Falls-out      & 0.594/0.595/0.596 & 0.595/0.595/0.595 \\
\hline
Regular 10-10-10     & 0.167/0.174/0.174 & 0.056/0.092/0.093 \\
\hline
Regular 10-10-10-out & 0.224/0.226/0.226 & 0.099/0.119/0.119 \\
\hline
Regular 10-10-30     & 0.161/0.164/0.165 & 0.132/0.134/0.134 \\
\hline
Regular 10-10-30-out & 0.155/0.161/0.164 & 0.256/0.258/0.259 \\
\hline
Chicago 10           & 0.241/0.249/0.249 & 0.355/0.363/0.363 \\
\hline
Chicago 10-out       & 0.256/0.264/0.264 & 0.365/0.373/0.373 \\
\hline
\end{tabular}
\label{tab:EMm}
\end{table}

% % -------------------------------------------------------
% % Mean arc flows
% % -------------------------------------------------------
% \begin{table}[h!]
% \centering
% \caption{Algorithm Performance -- Mean arc flows}
% \begin{tabular}{| l | c | c | c | c | c |}
% \hline
% \textbf{Network} 
%   & \textbf{Lin. Reg.} 
%   & \textbf{NN (MSE)} 
%   & \textbf{NN} 
%   & \textbf{NN+FP} 
%   & \textbf{M-P} \\
% \hline
% Small                & 12.355  & -- & -- & -- & 11.434 \\
% \hline
% Small-out            & 13.626  & -- & -- & -- & 12.089 \\
% \hline
% Steenbrink           & 21.987  & -- & -- & -- & 19.883 \\
% \hline
% Steenbrink-out       & 24.367  & -- & -- & -- & 21.785 \\
% \hline
% Nguyen Dupuis        & 27.916  & -- & -- & -- & 27.862 \\
% \hline
% Nguyen Dupuis-out    & 26.913  & -- & -- & -- & 26.913 \\
% \hline
% Sioux Falls          & 730.683 & -- & -- & -- & 14.422 \\
% \hline
% Sioux Falls-out      & 757.509 & -- & -- & -- & 15.837 \\
% \hline
% Regular 10-10-10     & 24.177  & -- & -- & -- & 24.176 \\
% \hline
% Regular 10-10-10-out & 25.765  & -- & -- & -- & 25.765 \\
% \hline
% Regular 10-10-30     & 66.963  & -- & -- & -- & 66.913 \\
% \hline
% Regular 10-10-30-out & 68.527  & -- & -- & -- & 68.527 \\
% \hline
% Chicago 10           & 1.089   & -- & -- & -- & 1.089  \\
% \hline
% Chicago 10-out       & 1.096   & -- & -- & -- & 1.095  \\
% \hline
% \end{tabular}
% \label{tab:Prom}
% \end{table}

% -------------------------------------------------------
% Inference Time
% -------------------------------------------------------
\begin{table}[h!]
\footnotesize
\centering
\caption{Average inference time (seconds)}
\begin{tabular}{| l | c | c | c | c | c | c | c | c |}
\hline
\multirow{2}{*}{\textbf{Network}} 
& \multirow{2}{*}{\textbf{Scilab}} 
  & \multicolumn{3}{c|}{\textbf{Linear Regression}} 
  & \multirow{2}{*}{\textbf{NN (MSE)}} 
  & \multirow{2}{*}{\textbf{NN}} 
  & \multirow{2}{*}{\textbf{NN+FP}} 
  & \multirow{2}{*}{\textbf{M-P}} \\
\cline{3-5}
& & \textbf{Classic} & \textbf{Ridge} & \textbf{Lasso} & & & & \\
\hline
Small           & 0.079712 & 0.000002 & 0.000001 & 0.000001 & 0.000002 & 0.000017 & 0.005000 & 0.00003 \\
\hline
Steenbrink      & 0.097809 & 0.000002 & 0.000002 & 0.000001 & 0.000003 & 0.000022 & 0.006000 & 0.00008 \\
\hline
Nguyen Dupuis   & 0.031645 & 0.000002 & 0.000001 & 0.000001 & 0.000003 & 0.000059 & 0.002800 & 0.00006 \\
\hline
Regular 10-10   & 0.098870 & 0.000001 & 0.000002 & 0.000001 & 0.000003 & 0.000003 & 0.007721 & 0.00018 \\
\hline
Regular 10-30   & 0.171582 & 0.000003 & 0.000001 & 0.000002 & 0.000003 & 0.000013 & 0.017570 & 0.00092 \\
\hline
Sioux Falls     & 0.106485 & 0.000003 & 0.000003 & 0.000003 & 0.000006 & 0.000685 & 0.017501 & 0.00098 \\
\hline
Chicago 10      & 0.057330 & 0.000120 & 0.000126 & 0.000144 & 0.000005 & 0.000631 & 0.875375 & 0.01208 \\
\hline
\end{tabular}
\label{tab:time}
\end{table}


\begin{thebibliography}{99}

\bibitem{Beckmann} Beckmann, M.J.,  On the Theory of Traffic Flows in Networks. Traffic Quarterly, 21, 109-116, 1967.

\bibitem{Bell-Iida} Bell M. G., Iida Y. Transportation Network Analysis, John Wiley and Sons, Chichester, U.K., 1997.

\bibitem{Choetal} Cho H.J., Smith T. E., Friesz T. L., A reduction method  for local sensitivity analyses of network equilibrium arc flows, Transportation Research, 34 B, pp. 31-51, 2000.

\bibitem{DafermosNagurney} Dafermos, S. and Nagurney, A., Sensitivity analysis for the asymmetric network equilibrium problem, Mathematical Programming, 28 (2), pp. 174-184, 1984.

\bibitem{DafermosSparrow} Dafermos, S. and Sparrow, F. T., Traffic assignment problem for a general network, Journal of Research of the National Bureau of Standards, 73 B, pp. 91-118, 1969.

\bibitem{Fan} Fan, W., Tang, Z., Ye, P., Xiao, F. and Zhang, J., Deep learning-based dynamic traffic assignment with incomplete origin–destination data. Transportation Research Record, 2677(3), pp. 1340–1356, (2023).

\bibitem{Graf} Graf, L., Harks, T., Kollias, K. ans Markl, M., Machine-learned prediction equilibrium for dynamic traffic assignment. Proceedings of AAAI, 36(5), pp. 5059–5067, 2022. 

\bibitem{WN} Jungel, K., Paccagnan, D., Parmentier, A. and Schiffer, M., WardropNet: Traffic flow predictions via equilibrium-augmented learning,  arXiv: 2410.06656, 2024.

\bibitem{Klein} Klein, I., Levy, N. and Ben-Elia, E., An agent-based model of the emergence of cooperation and a fair and stable system optimum using ATIS on a simple road network, Transportation Research Part C: Emerging Technologies, 86, pp. 183-201, 2018.

\bibitem{endtoend} Liu Z., Yin Y., Fan B., Grimm D. K., End-to-End Learning of User Equilibrium  with Implicit Neural Networks, Transportation Research Part C: Emerging Technologies, 150, 2023.

\bibitem{LiuGN} Liu, T. and Meidani, H.,  End-to-end heterogeneous graph neural networks for traffic assignment, Transportation Research Part C: Emerging Technologies, 2024.

\bibitem{endtoend25} Liu, Z. and Yin, Y., End-to-end learning of user equilibrium: Expressivity, generalization, and optimization, Transportation Science, 59(4), pp. 853–882, 2025. 

\bibitem{Heaton} McKenzie, D., Heaton, H., Li, Q., Wu, F., Osher, S. and Yin, W., Three-Operator Splitting for Learning to Predict Equilibria in Convex Games, SIAM Journal on Mathematics of Data Science, 6(3), pp.627-648, 2024.

\bibitem{Morandi} Morandi V., Bridging the user equilibrium and the system optimum in static traffic assignment: a review,  4OR-Q J Oper Res 22, pp.89–119, 2024. 


\bibitem{Outrata} Outrata J. V., On a special class of mathematical programs with equilibrium constraints. Lecture Notes Econom. Math Systems,  452., pp. 246-260, 1997.

\bibitem{PatRock} Patriksson M., Rockafellar R. T., Sensitivity analysis of aggregated variational inequality problems, with applications to traffic equilibria, Transportation Science, 37, pp. 56-68, 2003.

\bibitem{Pat04} Patriksson M., Sensitivity analysis of Traffic Equilibria, Transportation Science, 38, pp. 258-281, 2004.

\bibitem{Patriksson} Patriksson M., The Traffic Assignment Problem: Models and Methods,  Dover Publications, 2015.

\bibitem{Qiuetal} Qiu Y., Magnanti T. L., Sensitivity analysis for variational inequalities defined on polyhedral sets , Math. Oper. Research,  14, pp. 410-432, 1989.

\bibitem{RH} Rahman, R. and Hasan, S., Data-driven traffic assignment: A novel approach for learning traffic flow patterns using graph convolutional neural network. Data Science for Transportation, 5(2), 11, 2023.

\bibitem{Tobin-Friesz} Tobin R. L., Friesz T. L., Sensitivity analysis for equilibrium network flow, Transportation Science,  22, pp. 242-250, 1988.

\bibitem{Wegner} Wegner S. V., Mathematical introduction to Data Science, Springer, Hamburg, Germany, 2024.

\bibitem{Yen} Yen N. D., Lipschitz continuity of solutions of variational inequalities with a parametric  polyhedral constraint,  Math. Oper. Research,  20, 695-708.

\end{thebibliography}
\end{document}